\documentclass[11pt]{amsart}
\usepackage[T1]{fontenc}
\usepackage{graphicx} 
\usepackage{amsmath}
\usepackage{amsthm} 
\usepackage{amssymb}
\usepackage{verbatim} 
\usepackage{color}
\usepackage[colorlinks=true, linkcolor=purple, urlcolor=purple, citecolor=purple]{hyperref}
\usepackage{enumitem}
\usepackage{tikz}
\usepackage{tikz-cd}
\usetikzlibrary{positioning}
\usepackage{forloop}

\newtheorem{theorem}{Theorem}[section]

\newtheorem{lemma}[theorem]{Lemma}

\theoremstyle{definition}
\newtheorem{definition}[theorem]{Definition}

\theoremstyle{remark}
\newtheorem{claim}{Claim}

\title{Equivalence of Cohen generics is hyperfinite}
\author{Iian B.~Smythe}
\address{Department of Mathematics and Statistics, University of Winnipeg, 515 Portage Avenue, Winnipeg, Manitoba, Canada, R3B 2E9}
\email{i.smythe@uwinnipeg.ca}
\thanks{The author is partially supported by an NSERC Discovery Grant (RGPIN-2023-03343).}

\date{\today}

\begin{document}

\begin{abstract}
    It is shown that the equivalence relation on Cohen generic reals over a countable model of set theory, where two are equivalent if they produce the same generic extension, is hyperfinite, answering a question of the author. The proof uses asymptotic separation and is a digestion of an argument originally produced by GPT-6.
\end{abstract}

\maketitle

\section{Introduction}

Throughout, let $M$ be a countable transitive model of $\mathsf{ZFC}$, or whichever finite fragment is sufficient to carry out all of the arguments herein. In \cite{MR4130394}, Clemens, Coskey and Dworetzky, as part of their study of the complexity of the isomorphism relation on countable models of set theory, considered the equivalence relation $\equiv^\mathbb{C}_M$ on Cohen generic reals over $M$ given by:
\[
    x\equiv^\mathbb{C}_M y \quad\Longleftrightarrow\quad M[x]=M[y].
\]
They showed, in their proof of \cite[Theorem 3.2]{MR4130394}, that $\equiv^\mathbb{C}_M$ is non-smooth, that is, not Borel reducible to the equality relation on any Polish space.

This example was expanded upon by the author in \cite{MR4452126}, where he developed the following general framework: Given a notion of forcing $\mathbb{P}$ in $M$ (meaning $\mathbb{P}$ and its ordering $\leq$ are both elements of $M$), there is a Polish subspace $\mathrm{Gen}^\mathbb{P}_M$ of $2^\mathbb{P}$ consisting of all $M$-generic filters for $\mathbb{P}$, on which the equivalence relation given by:
\[
    G\equiv^\mathbb{P}_M H \quad\Longleftrightarrow\quad M[G]=M[H].
\]
is a countable Borel equivalence relation. In the case of Cohen forcing $\mathbb{C}$, the author showed that $\equiv^\mathbb{C}_M$ is hyperhyperfinite \cite[Theorem 3.1]{MR4452126}, and asked \cite[Question 2]{MR4452126} whether it is hyperfinite. Recall that an equivalence relation is \emph{hyperfinite} if it is a countable increasing union of \emph{finite} (i.e., having finite classes) Borel equivalence relations, and is \emph{hyperhyperfinite} if it is a countable increasing union of hyperfinite equivalence relations. 

It it is a long-standing open problem, going back at least to \cite{MR1149121}, known as ``the union problem'', whether every hyperhyperfinite equivalence relation is hyperfinite. The question of whether the equivalence of Cohen generics $\equiv^\mathbb{C}_M$ is a possible counterexample to the union problem has been raised by several authors, for example in \cite{barmpalias2026openproblemsmathematicallogic}, \cite{calderoni2026forcinggenericitycbers}, and \cite{MR4837611}. In this note, it is shown that this is not the case:

\begin{theorem}\label{thm:cohen_is_hf}
    $\equiv^\mathbb{C}_M$ is hyperfinite.
\end{theorem}

The proof of Theorem \ref{thm:cohen_is_hf} is an application of the theory of Borel asymptotic dimension developed by Conley, Jackson, Marks, Seward, and Tucker-Drob in \cite{MR4679959}, together with ideas from the original hyperhyperfiniteness argument in \cite{MR4452126}.

\subsection*{AI disclosure:} The proof of Theorem \ref{thm:cohen_is_hf} was originally produced by GPT-6 Astra in response to prompts from Zolt\'{a}n Vidny\'{a}nszky instructing it to prove that $\equiv^\mathbb{C}_M$ is hyperfinite using an asymptotic dimension argument. After GPT produced a solution, Vidny\'{a}nszky communicated it to the author, who independently verified its correctness. It was agreed by both that the AI-generated solution should be translated into a more readable form and communicated as soon as possible. This note is the result, and was written entirely by the author, who takes responsibility for its content and correctness in accordance with the Leiden Declaration \cite{LeidenDeclaration2026}.

\section{Automorphisms and Borel codes}

If $\mathbb{P}$ is a complete Boolean algebra in $M$,\footnote{Meaning $M\models\text{``$\mathbb{P}$ is a complete Boolean algebra''}$.} then it follows by a result of Vop\v{e}nka and H\'{a}jek \cite{MR444473} (see also \cite[Theorem 3.5.1]{MR373889}) that $\equiv^\mathbb{P}_M$ is induced by the natural action of the countable group $\mathrm{Aut}^M(\mathbb{P})=\mathrm{Aut}(\mathbb{P})\cap M$ of automorphisms of $\mathbb{P}$, as computed in $M$, on the space of its $M$-generic filters. Since every poset densely embeds into a complete Boolean algebra, and dense embeddings yield Borel isomorphic equivalence relations \cite[Lemmas 2.7]{MR4452126}, one may, without loss of generality, work with complete Boolean algebras.

As in \cite{MR4452126}, when analyzing the complexity of $\equiv^\mathbb{C}_M$, it will be convenient to use a particular presentation of Cohen forcing and its Boolean completion. We first need to recall the notion of a \emph{Borel code}, originally developed in \cite{MR265151}. If $M$ thinks $B$ is a Borel subset of $2^\omega$, say, then $M$ contains a real, the \emph{code} for $B$, which describes how $B$ is built up from a fixed computable basis of open sets. We write $B^M$ for the interpretation of that code in $M$, and omit the superscript and just write $B$ for the corresponding interpretation in the ambient universe $\mathbf{V}$. We will rely on the absoluteness of various properties of Borel sets (or functions) coded in $M$, meaning the property holds of $B^M$ in $M$ if and only if it holds of $B$ in $\mathbf{V}$.

\begin{lemma}\label{lem:absolute1}
    For Borel sets, functions, or equivalence relations coded in $M$, the following notions are absolute:
    \begin{enumerate}
        \item For sets: membership, containment, countable unions and intersections, meagerness;
        \item For functions: being a function, images, injectivity, surjectivity;
        \item For equivalence relations: being an equivalence relation, having finite classes, being induced by the orbit equivalence relation of the group generated by functions $f_0,\ldots,f_{n-1}$ coded in $M$, having a function $s$ as a selector.
    \end{enumerate}
\end{lemma}

\begin{proof}
    All of these statements, with the exception of the last two claims in (3), are proved in \cite{MR4452126} or \cite{MR265151} by showing that the associated property admits a $\Pi^1_1$ definition in the appropriate codes. The second-to-last claim in (3) is similar to \cite[Lemma 3.5(iii)]{MR4452126}. For the remaining claim, a function $s$ is a selector for a Borel equivalence relation $E$ if and only if $$\forall x,y [xEs(x)\land (xEy\rightarrow s(x)=s(y))],$$ which is clearly $\Pi^1_1$ in the appropriate codes, thus absolute.
\end{proof}

In $M$, let $\mathbb{C}$ be the collection of all non-meager Borel subsets of $2^\omega$, and $\overline{\mathbb{C}}$ the collection of Borel sets in $2^\omega$, modulo the ideal of meager sets. Then, $\overline{\mathbb{C}}$ is a complete Boolean algebra in $M$ and the natural map $i:\mathbb{C}\to\overline{\mathbb{C}}$, mapping a set to its equivalence class modulo meager, is a dense embedding. If $G\subseteq\mathbb{C}$ is a filter for $\mathbb{C}$, we write $\overline{G}$ for the filter generated by $i''G$ in $\overline{\mathbb{C}}$. Then, $G$ is $M$-generic if and only if $\overline{G}$ is, and $M[G]=M[\overline{G}]$. (See \cite[\S VII.7]{MR597342}.)

Back in $\mathbf{V}$, we call $x\in 2^\omega$ an \emph{$M$-generic Cohen real} if there is an $M$-generic filter $G$ for $\mathbb{C}$ such that
\[
    \{x\}=\bigcap\{B:B^M\in G\},
\]
in which case $M[x]=M[G]$. Equivalently, $x$ is in every comeager Borel subset of $2^\omega$ which is coded in $M$. Since $G$ is uniquely definable from $x$ as
\[
    G=\{B^M:x\in B\},
\]
we may write $G=G_x$. Let $\mathcal{C}_M$ be the set of all $M$-generic Cohen reals and write $\equiv^\mathbb{C}_M$ for relation
\[
    x\equiv^\mathbb{C}_M y \quad\Longleftrightarrow\quad M[x]=M[y]
\]
on $\mathcal{C}_M$. Then, by \cite[Lemma 2.10]{MR4452126}, the mapping $x\mapsto G_x$ is a Borel isomorphism between our two definitions of $\equiv^\mathbb{C}_M$, one on $M$-generic Cohen reals, the other on $M$-generic filters for $\mathbb{C}$.

A result of Maharam and Stone \cite[Corollary 2]{MR522175} says that every automorphism of the algbera of Borel subsets of $2^\omega$ modulo meager sets is induced by a meager-class preserving (meaning both $g$ and $g^{-1}$ preserve meager sets) Borel bijection $g:2^\omega\to 2^\omega$. Let $\mathrm{MCP}^M(2^\omega)$ be the set of all meager-class preserving Borel bijections $2^\omega\to 2^\omega$ which are coded in $M$. Clearly, $\mathrm{MCP}^M(2^\omega)$ forms a countable group and, by the absoluteness of meagerness, acts on $\mathcal{C}_M$ in a way which preserves $\equiv^\mathbb{C}_M$. It likewise acts on the space of $M$-generic filters for $\mathbb{C}$ by
\[
    g\cdot G=\{(g^M)''(B):B\in G\},
\]
for $g\in\mathrm{MCP}^M(2^\omega)$ and $G\in\mathrm{Gen}^\mathbb{C}_M$. The following lemma is implicit in the proof of hyperhyperfiniteness of $\equiv^\mathbb{C}_M$ in \cite{MR4452126}; we make it explicit here:

\begin{lemma}\label{lem:realize}
    $\equiv^\mathbb{C}_M$ is induced by the action of $\mathrm{MCP}^M(2^\omega)$ on $\mathcal{C}_M$.
\end{lemma}

\begin{proof}
    It will suffice to show that, for $x,y\in\mathcal{C}_M$, $M[x]=M[y]$ if and only if there $g\in\mathrm{MCP}^M(2^\omega)$ such that $g(x)=y$. We proceed with a series of equivalences:
    \begin{align*}
        M[x]=M[y]&\quad\Longleftrightarrow\quad M[\overline{G_x}]=M[\overline{G_y}]\\
        &\quad\Longleftrightarrow\quad \exists\gamma\in\mathrm{Aut}^M(\overline{\mathbb{C}})(\gamma''\overline{G_x}=\overline{G_y})\\
        &\quad\Longleftrightarrow\quad \exists g\in\mathrm{MCP}^M(2^\omega)(g\cdot G_x=G_y)\\
        &\quad\Longleftrightarrow\quad \exists g\in\mathrm{MCP}^M(2^\omega)(G_{g(x)}=G_y)\\
        &\quad\Longleftrightarrow\quad \exists g\in\mathrm{MCP}^M(2^\omega)(g(x)=y).
    \end{align*}
    The first equivalence follows from the fact that $M[x]=M[G_x]=M[\overline{G_x}]$ and $M[y]=M[G_y]=M[\overline{G_y}]$. The second uses the completeness (in $M$) of $\overline{\mathbb{C}}$ and the Vop\v{e}nka--H\'{a}jek  theorem, while the third follows by the Maharam--Stone theorem. The fourth equivalence is by \cite[Lemma 3.7]{MR4452126}, and the last by the fact that an $M$-generic filter for $\mathbb{C}$ uniquely determines the corresponding Cohen real.
\end{proof}

\section{Generic asymptotic separation}

If $H$ is a locally finite Borel graph on a standard Borel space $X$, denote by $\rho_H$ the extended path metric on $X$, where $\rho_H(x,y)=\infty$ if $x$ and $y$ are not in the same connected component of $H$. For any integer $r\geq 1$, let $H^{\leq r}$ be the graph which joins all points of $X$ which are at most $\rho_H$-distance $r$ apart. In \cite{MR4679959}, Conley, Jackson, Marks, Seward, and Tucker-Drob introduced \emph{Borel asymptotic dimension}, a Borel form of the notion of asymptotic dimension, and \emph{asymptotic separation index}, for Borel extended metric spaces. It is the latter which will concern us here.

\begin{definition}
    Let $H$ be a locally finite Borel graph on a standard Borel space $X$. The \emph{asymptotic separation index} of $H$, written $\mathrm{asi}(H)$, is the least $d\in\omega$ such that for every $r\geq 1$, there are Borel sets $U_0,U_1,\ldots, U_d$ which cover $X$ and such that the connected components of $H^{\leq r}\upharpoonright U_i$ are finite, for all $i\leq d$.
\end{definition}

The following result, which we state in the special case of a locally finite Borel graph (see \cite[Theorem 1.15]{MR4969378} for this formulation), says that locally finite Borel graphs have asymptotic separation index $\leq 1$, generically.

\begin{theorem}[{\cite[Theorem 4.8(b)]{MR4679959}}]\label{thm:asi}
    If $H$ is a locally finite Borel graph on a standard Borel space $X$, then there is an $H$-invariant comeager Borel set $D\subseteq X$ such that $\mathrm{asi}(H\upharpoonright D)\leq 1$.
\end{theorem}

We note here the absoluteness of the relevant notions for Borel graphs.

\begin{lemma}\label{lem:absolute2}
    For $H$ a Borel set coded in $M$, the following are upwards absolute (meaning if they hold of the coded set in $M$, then they hold in $\mathbf{V}$):
    \begin{enumerate}
        \item $H$ is a locally finite Borel graph;
        \item $\rho_H$ is the extended path metric on $H$;
        \item for $r\geq 1$, the graph $H^{\leq r}$;
        \item the connectedness relation $R$ on $H$.
    \end{enumerate}
\end{lemma}

\begin{proof}
    These are proved in much the same way as Lemma \ref{lem:absolute1}. The only wrinkle is that whether $x$ and $y$ are in the same connected component of $H$ requires asserting the existence of a path (i.e., finitely many points) from $x$ to $y$; this existential quantifier makes the connectedness relation (explicit in (4), but implicit in the definition of $\rho_H$, and thus in (2) and (3)) upwards absolute, which is all we require.
\end{proof}

\section{Proof of hyperfiniteness}

The following is the main result of this note.

\begin{theorem}\label{thm:general_hf}
    Let $\Gamma=\{f_n:n\in\omega\}$ be a countable group of meager-class preserving Borel bijections on $2^\omega$ such that each $f_n$ is coded in $M$. Then, the orbit equivalence relation induced by $\Gamma$ on $2^\omega$ is hyperfinite when restricted to the space of $M$-generic Cohen reals $\mathcal{C}_M$.
\end{theorem}

Theorem \ref{thm:cohen_is_hf} is an immediate consequence of Theorem \ref{thm:general_hf}, Lemma \ref{lem:realize}, and the fact that $\mathrm{MCP}^M(2^\omega)$ is countable. Note that while we require each $f_n$ to be coded in $M$, the whole sequence $(f_n)_{n\in\omega}$ need not be, as for instance, $\mathrm{MCP}^M(2^\omega)$ is uncountable in $M$. 

The proof of Theorem \ref{thm:general_hf}, as produced by GPT, is a combination of the proof that $\equiv^\mathbb{C}_M$ is hyperhyperfinite from \cite{MR4452126}, with Theorem \ref{thm:asi} now playing the role of the generic hyperfiniteness theorem  \cite[Theorem 12.1]{MR2095154}, together with the proof that locally finite graphs of finite Borel asymptotic dimension induce hyperfinite Borel equivalence relations \cite[Theorem 7.1]{MR4130394}.

\begin{proof}[Proof of Theorem \ref{thm:general_hf}]
    Let $\Gamma$ be as described and $E$ the orbit equivalence relation induced by $\Gamma$, restricted to the set $\mathcal{C}_M$.
    
    In $M$ (we will omit the superscripts while working in $M$ here, for notational brevity), let $R_n$ be the orbit equivalence relation on $2^\omega$ induced by the group generated by $\{f_0,\ldots,f_{n-1}\}$. Then, the $R_n$'s form an increasing sequence of countable Borel equivalence relations in $M$, and in $\mathbf{V}$,
    \[
        E=\bigcup_{n\in\omega}R_n\upharpoonright\mathcal{C}_M.  
    \]
    
    We will construct an increasing sequence of finite Borel equivalence relations $F_n$ on $2^\omega$ in such a way that each is coded in $M$ and $F_n\subseteq R_n$. We begin by letting $F_0=\Delta_{2^\omega}$, the equality relation on $2^\omega$.
    
    For $n\geq 1$, suppose that $F_{n-1}$ has been defined as required. We now work in $M$. Choose a Borel selector $s_{n-1}$ for $F_{n-1}$. Let
    \[
        H_n=\left(F_{n-1}\cup\bigcup_{i=0}^{n-1}(\mathrm{graph}(f_i)\cup\mathrm{graph}(f_i^{-1}))\right)\setminus\Delta_{2^\omega}.
    \]
    Then, $H_n$ is a locally finite Borel graph on $2^\omega$ in which each $F_n$-class is a clique, and its connectedness relation coincides with $R_n$. 
    
    We apply Theorem \ref{thm:asi} in $M$ to obtain Borel sets $D_n, A_n^0, A_n^1\subseteq 2^\omega$ such that $D_n$ is comeager in $2^\omega$, $D=A_n^0\sqcup A_n^1$, and every connected component of $H_n^{\leq 4}\upharpoonright A_n^i$ is finite, for $i=0,1$. Let $P_n$ be the finite Borel equivalence relation on $2^\omega$ where, on each $A_n^i$, its classes are the connected components of $H_n^{\leq 4}$, and outside $D_n$ it is equality. Note that $P_n\subseteq R_n$.
    
    Define $F_n$ on $2^\omega$ by:
    \[
        xF_n y \quad\iff\quad s_{n-1}(x) P_n s_{n-1}(y).
    \]  
    Then, $F_n$ is a finite Borel equivalence relation with $F_{n-1}\subseteq F_n$, since it is the pullback of a finite Borel equivalence relation under the selector $s_{n-1}$ for $F_{n-1}$, and $F_n\subseteq R_n$. To see the latter, note that if $xF_n y$, then
    \[
        xF_{n-1} s_{n-1}(x),\quad s_{n-1}(x)P_ns_{n-1}(y),\quad\text{and}\quad s_{n-1}(y)F_{n-1}y,
    \]
    and since both $F_{n-1}$ and $P_n$ are subsets of $R_n$, we have that $xR_n y$.

    This describes the construction of the $F_n$'s. We now work in $\mathbf{V}$, using Lemmas \ref{lem:absolute1} and \ref{lem:absolute2} for absoluteness of all relevant notions. Since all of the $D_n$'s are comeager Borel sets coded in $M$, we have that $\mathcal{C}_M\subseteq\bigcap_{n\in\omega}D_n$.

    \begin{claim}
        For each $x\in\mathcal{C}_M$ and $n\in\omega$, the partial orbit 
        \[
            T_n=\{x,f_0(x),\ldots,f_{n-1}(x)\}
        \]
        meets at most two $F_n$-classes.
    \end{claim}

    If $u\in T_n(x)$, then $\rho_{H_n}(x,u)\leq 1$, and since $u$ and $s_{n-1}(u)$ lies in the same $F_{n-1}$-class, so $\rho_{H_n}(s_{n-1}(u),u)\leq 1$ as well. Thus, for any $u,v\in T_n(x)$, 
    \begin{align*}
        \rho_{H_n}(s_{n-1}(u),s_{n-1}(v))\leq\;&\rho_{H_n}(s_{n-1}(u),u)+\rho_{H_n}(u,x)\\
        &+\rho_{H_n}(x,v)+\rho_{H_n}(v,s_{n-1}(v))\leq 4. 
    \end{align*}
    Note that $s_{n-1}(u),s_{n-1}(v)\in\mathcal{C}_M$ since they lie in the same $R_n$-class as $x$, and so are in $D_n=A_n^0\sqcup A_n^1$. If both $s_{n-1}(u),s_{n-1}(v)$ are in the same $A_n^i$, then they are equal or joined by an edge of $H_n^{\leq 4}$, in which case $s_{n-1}(u)P_n s_{n-1}(v)$, and thus $uF_n v$. It follows that $T_n$ meets at most two $F_n$-classes, as claimed.

    Define an equivalence relation on $\mathcal{C}_M$ 
    \[
        F=\bigcup_{n\in\omega}F_n\upharpoonright\mathcal{C}_M.
    \]
    By construction, $F$ is hyperfinite and $F\subseteq \bigcup_{n\in\omega}R_n\upharpoonright\mathcal{C}_M=E$.

    \begin{claim}
        Each $E$-class meets at most two $F$-classes.
    \end{claim}

    Suppose that $y_0,y_1,y_2$ are members of the $E$-class of some $x$. Choose $n$ large enough so that $y_0,y_1,y_2\in T_n(x)$. Then, by Claim 1, two of them are in the same $F_n$-class, and hence the same $F$-class. Thus, an $E$-class cannot contain elements from more than two distinct $F$-classes.

    In other words, $F$ is an index $\leq 2$ hyperfinite subequivalence relation of $E$, and thus, by \cite[Proposition 1.3(vii)]{MR1900547}, $E$ is hyperfinite.
\end{proof}

\newpage

\bibliographystyle{abbrv}
\bibliography{cohen.bib}

\end{document}